\documentclass[12pt,leqno,twoside]{article}
\usepackage{dh,latexsym}

\title{\bf Bounds on norms of compound matrices and on products
of eigenvalues \thanks{The research
was started while the first author visited the Technion in 1995.
It was continued during a visit of the second author at
Universit\"{a}t
Bielefeld and at the University of Wisconsin - Madison in 1997,
and it was completed during a visit of the third author at the
Technion in 1997. The research was partly supported by SFB 343
"Diskrete Strukturen in der Mathematik", Bielefeld, and by NSF
Grant
DMS-9424346.}}

\author{Ludwig Elsner\\
\it Fakult\"{a}t f\"{u}r Mathematik\\
\it Universit\"{a}t Bielefeld\\
\it Postfach 100131\\
\it D-33615 Bielefeld\\
\it Germany
\and
Daniel Hershkowitz\\
\it Department of Mathematics\\
\it Technion - Israel Institute of Technology\\
\it Haifa 32000\\
\it Israel
\and
Hans Schneider\\
\it Department of Mathematics\\
\it University of Wisconsin\\
\it Madison, Wisconsin 53706\\
\it U.S.A}

\begin{document}

\bibliographystyle{alpha}

\date{21 November 1997}

\maketitle

\begin{abstract}

An upper bound on operator norms of compound matrices is
presented, and
special cases that involve the $\ell_1$, $\ell_2$ and
$\ell_\infty$ norms are investigated. The results are then used
to obtain bounds on products of the largest or smallest
eigenvalues of a matrix.

\end{abstract}

\section{Introduction \label{intro}}

Let $A$ be a complex matrix and let $C_k(A)$ be its $k$th
compound. It was
shown in \cite[Formula (12)]{Sch} that the maximal row sum (of
moduli)
of elements of $C_k(A)$ is less than or equal to the product of
the $k$ largest rows sums of $A$, and it follows that the product
of $k$ largest (moduli of) eigenvalues of $A$ is bounded above by
the product of the $k$
largest row sums of $A$. The case of equality in these
inequalities
investigated in \cite[Theorems I and II]{O}.

\bigskip

The results in \cite{Sch} and \cite{O} can be viewed as relating
the $\ell_1$ norm of rows of a matrix to the $\ell_1$ norm of its
compounds (viewed as an operator on rows). Working in terms of
columns, we consider in this paper the relations between other
norms $\nu$ of columns and norms $\mu$ of the compounds. We begin
by proving a general result of the above type which involves a
constant $\theta_k(\mu,\nu)$.  We evaluate this constant in some
special cases that involve the $\ell_1$, $\ell_2$ and
$\ell_\infty$ norms. Again , this leads naturally to upper bounds
on the product of the $k$ largest eigenvalues or, equivalently,
lower bounds on the product of the $k$ smallest eigenvalues,
which involve products of norms of columns and of norms of rows
of the matrix.

\bigskip

As a consequence of our theorems, we obtain generalizations of
results of \cite{HP}, \cite{SW} and \cite{Smi} on bounds on norms
of the adjoint matrix, which is essentially the $n-1$ compound
matrix, to $k$th compound matrices. The application of our
theorems to the adjoint case sharpens the results in \cite{HP},
\cite{SW} and \cite{Smi}.

\section{Upper bounds on norms of compound matrices
\label{upper}}

Let $A$ be a matrix in $\C^{nn}$.
For subsets $\alpha$ and $\beta$ of $\{1,\ldots,n\}$
we denote by $A(\alpha|\beta)$ the submatrix of $A$ whose rows
are
indexed by $\alpha$ and whose columns are indexed by $\beta$ in
their
natural order.

\bigskip

Let $k$ be a positive integer, $k \leq n$. We denote by
$C_k(A)$ the $k$th compound of the matrix $A$, that is, the
${n \choose k} \times {n \choose k}$ matrix whose elements are
the minors $\det A(\alpha|\beta)$, $\alpha,\beta \subseteq
\{1,\ldots,n\}$,
$|\alpha|=|\beta|=k$. We index $C_k(A)$ by $\alpha \subseteq
\{1,\ldots,n\}$, $|\alpha|=k$ (ordered lexicographically).

\bigskip

Let $\nu$ be a vector norm on $\C^n$, and for a positive integer
$k,\,\,k \leq n$, let $\mu$ be a (submultiplicative) norm on
$\C^{mm}$ where $m={n \choose k}$. We define
$$
\theta_k (\mu,\nu) = \max \{\mu(C_k(B)): B \in
\C^{nn},\;\nu(col_i(B))=1,\;i=1,\ldots,n\},
$$
where $col_i(B)$ denotes the $i$th column of $B$.

\bigskip

The following theorem is the main tool from which we derive our
results.

\begin{thm}
\label{thm0}
For an absolute operator norm $\mu$ we have
\beq
\label{a0}
\mu(C_k(A)) \leq \theta_k(\mu,\nu)
\max_{\stackrel{\scriptstyle \alpha \subseteq
\{1,\ldots,n\}}{|\alpha|=k}} \prod_{i \in \alpha}\nu(col_i(A)).
\end{equation}
\end{thm}

{\it Proof.} Assume first that $A$ has no zero columns.
Note that for every nonsingular matrix $R$ we have
$$
C_k(AR^{-1})C_k(R) = C_k(A).
$$
In particular, if we choose the matrix $R$ to be the diagonal
matrix
$diag(r_i)$ where $r_i = \nu(col_i(A)),\,\,i=1,\ldots,n$,
then $C_k(R)$ is a diagonal matrix with diagonal elements
$\prod_{i \in \alpha}r_i$, $\alpha \subseteq \{1,\ldots,n\}$,
${|\alpha|=k}$. Since $\mu$ is an absolute operator norm, it
follows by Theorem 3 in \cite{BSW}, see also \cite[p.310, Theorem
5.6.37]{HJ}, that
$$
\mu(C_k(R)) = \max_{\stackrel{\scriptstyle \alpha \subseteq
\{1,\ldots,n\}}{|\alpha|=k}} \prod_{i \in \alpha}r_i.
$$
Since the $\nu$ norm of the columns of $AR^{-1}$ are equal to
$1$, it now follows that
$$
\mu(C_k(A)) \leq \mu(C_k(AR^{-1}))\,\mu(C_k(R)) \leq
\theta_k(\mu,\nu)
\max_{\stackrel{\scriptstyle \alpha \subseteq
\{1,\ldots,n\}}{|\alpha|=k}} \prod_{i \in \alpha}\nu(col_i(A)).
$$
If $A$ has zero columns then we apply the above procedure to the
matrix $A + \epsilon I$ for $\epsilon$ sufficiently small, and
then use continuity arguments to prove our assertion.\eop

\begin{rem}
\rm Note that inequality (\ref{a0}) is sharp since, by definition
of $\theta_k(\mu,\nu)$, equality is attained for matrices whose
columns have $\nu$ norm $1$.
\end{rem}

\bigskip

In the rest of this section we evaluate $\theta_k(\mu,\nu)$ for
cases that involve the $\ell_p$ norms, $p=1,2,\infty$. We denote
these by $||.||_p$, and use the same notation for both vector
norms and matrix operator norms.

\bigskip

The following proposition is proven in \cite[Formulas (9) and
(10)]{O}. The inequality (\ref{a1}) also occurs in the proof of
Theorem 1 in \cite{Sch}.

\begin{pro}
\label{pro1}
For every $\beta \in \{1,\ldots,n\},\,\,|\beta|=k$ we have
\beq
\label{a1}
||col_\beta(C_k(A)||_1 \leq
\max_{\stackrel{\scriptstyle \alpha \subseteq
\{1,\ldots,n\}}{|\alpha|=k}} \prod_{i \in \alpha}||col_i(A)||_1.
\end{equation}
Furthermore, if $A$ has at least $k$ nonzero columns then for a
set
$\beta \in \{1,\ldots,n\},\,\,|\beta|=k$,
the following are equivalent:\\
{\rm (i)} Equality holds in {\rm (\ref{a1})}.\\
{\rm (ii)} We have
$$
\prod_{i \in \beta}||col_i(A)||_1 =
\max_{\stackrel{\scriptstyle \alpha \subseteq
\{1,\ldots,n\}}{|\alpha|=k}} \prod_{i \in \alpha}||col_i(A)||_1,
$$
and the columns of $A$ indexed by $\beta$ have disjoint
supports.
\end{pro}

\begin{defi}
\rm A matrix $A$ is said to be a {\it monomial} matrix if $A =
PD$, where $P$ is a permutation matrix and $D$ is a diagonal
matrix.
\end{defi}

\begin{thm}
\label{thml1}
We have $\theta_k(\ell_1,\ell_1) = 1.$
\end{thm}

{\it Proof.} Since
$$
||C_k(A)||_1 = \max_{\stackrel{\scriptstyle \beta \subseteq
\{1,\ldots,n\}}{|\beta|=k}}||col_\beta (C_k(A))||_1,
$$
it follows immediately from Proposition \ref{pro1} that
$\theta_k(\ell_1,\ell_1) \leq 1$.
By Proposition \ref{pro1}, every monomial matrix $A$ satisfies
$$
||C_k(A)||_1 = \max_{\stackrel{\scriptstyle \beta \subseteq
\{1,\ldots,n\}}{|\beta|=k}}||col_\beta (C_k(A))||_1 =
\max_{\stackrel{\scriptstyle \alpha \subseteq
\{1,\ldots,n\}}{|\alpha|=k}} \prod_{i \in \alpha}||col_i(A)||_1.
$$
implying that $\theta_k(\ell_1,\ell_1) = 1$.\eop

\begin{thm}
\label{tk2}
We have $\theta_k(\ell_2,\ell_2) =
\left(\frac{n}{k}\right)^\frac{k}{2}.$
Furthemore, if $A$ is nonsingular and $k<n$ then
\beq
\label{add2}
||C_k(A)||_2 < \left(\frac{n}{k}\right)^\frac{k}{2}
\max_{\stackrel{\scriptstyle \alpha \subseteq
\{1,\ldots,n\}}{|\alpha|=k}} \prod_{i \in \alpha}||col_i(A)||_2.
\end{equation}

\end{thm}

{\it Proof.} Let $A \in \C^{nn}$ be such that
$||col_i(A)||_2=1,\,\,\,i=1,\ldots,n$.
The matrix $B=A^*A$ is positive semidefinite with diagonal
entries equal to 1. Let $\rho_1\geq\rho_2\geq\ldots\rho_n\geq 0$
be the eigenvalues
of $B$. Note that $\rho_1+\ldots+\rho_n = trace(B) = n$. It now
follows that
\beq
\label{add1}
\left(||C_k(A)||_2\right)^2 = ||C_k(B)||_2 =
\rho_1\rho_2\cdots\rho_k \leq
\left(\frac{\rho_1+\ldots+\rho_k}{k}\right)^k
\leq\left(\frac{n}{k}\right)^k.
\end{equation}
We thus have
\beq
\label{k2}
\theta_k(\ell_2,\ell_2) \leq
\left(\frac{n}{k}\right)^\frac{k}{2}.
\end{equation}

\noindent
We now prove that equality holds in (\ref{k2}).
By \cite{Horn}, see also \cite[Theorem 2]{Mirsky}, there exists a
positive semidefinite
$n \times n$ matrix $B$ with diagonal elements all equal to $1$
and where
the eigenvalues of $B$ are $\frac{n}{k}$ with
multiplicity $k$ and $0$ with multiplicity $n-k$. We have
$$
||C_k(B)||_2 = \left(\frac{n}{k}\right)^k.
$$
Now, let $A$ be the positive semidefinite matrix such that $B =
A^2$. Since the diagonal entries of $B$ are all equal to $1$, it
follows that
$||col_i(A)||_2=1,\,\,\,i=1,\ldots,n$. Also,
$$
||C_k(A)||_2 =
\sqrt{||C_k(B)||_2} = \left(\frac{n}{k}\right)^\frac{k}{2},
$$
proving that equality holds in (\ref{k2}). Finally, notice that
if $A$ is nonsingular then $\rho_n > 0$ and so strict inequality
holds in (\ref{add1}) whenever $k<n$. Therefore, equality in
(\ref{k2}) cannot be attained for nonsingular matrices, and using
the techniques of the proof of Theorem \ref{thm0} one can prove
the strict inequality (\ref{add2}) whenever $A$ is nonsingular
and $k<n$.\eop

\bigskip

We remark that the inequality (\ref{k2}) in the case $k=n-1$ is
proven in \cite[Theorem 4]{Smi}, see also \cite[Lemma 2]{HP}. Our
proof of this inequality is essentially the same as in \cite{Smi}
and \cite{HP}. The equality case is, however, not handled in
these two references.

\begin{thm}
\label{thm212}
For $k<n$ we have
\beq
\label{a11}
\theta_k(\ell_\infty,\ell_\infty) \leq {n \choose k}
(k+1)^{\frac{k-1}{2}}.
\end{equation}
\end{thm}

{\it Proof.} Let $A \in \C^{nn}$ be such that
$||col_i(A)||_\infty=1,\,\,\,i=1,\ldots,n$.
Let $x$ be a vector in $\C^{n \choose k}$ and let $y = C_k(A)x$.
For
every subset $\alpha$ of $\{1,\ldots\,n\}$ of cardinality $k$ we
have
\beq
\label{e1}
y_\alpha = \sum_
{\stackrel{\scriptstyle \beta \subseteq
\{1,\ldots,n\}}{|\beta|=k}}
{\det A(\alpha|\beta)}x_\beta.
\end{equation}
We have $k<n$. Therefore, note that each subset $\beta$ of
$\{1,\ldots\,n\}$ of cardinality $k$
is contained in $n-k$ different subsets $\gamma$
of $\{1,\ldots\,n\}$ of cardinality $k+1$. Therefore, we have
\beq
\label{e2}
\sum_
{\stackrel{\scriptstyle \beta \subseteq
\{1,\ldots,n\}}{|\beta|=k}}
{\det A(\alpha|\beta)}x_\beta
= \frac{1}{n-k} \sum_
{\stackrel{\scriptstyle \gamma \subseteq
\{1,\ldots,n\}}{|\gamma|=k+1}}
\sum_
{\stackrel{\scriptstyle \beta \subseteq \gamma}{|\beta|=k}}
{\det A(\alpha|\beta)}x_\beta.
\end{equation}
Observe that the rightmost sum of (\ref{e2}) is the determinant
of the $(k+1) \times (k+1)$ matrix $B$ obtained by appending the
subvector of $x$ (with possible different
signs of elements) indexed by the subsets
$\beta$ of $\gamma$ of cardinality $k$ as a row to the matrix
$A(\alpha|\gamma)$.
Thus, if we choose $x$ such that $||x||_\infty = 1$ then the
matrix $B$
has entries of modulus less than or equal to $1$, and by the
Hadamard determinant theorem, e.g. \cite[p.114, Theorem
4.1.7]{MM} it follows that
$\det B \leq \sqrt{k+1}\,^{k+1}$.
Hence, it follows from (\ref{e1}) and (\ref{e2}) that
$$
||C_k(A)||_\infty \leq \frac{{n \choose
{k+1}}\sqrt{k+1}\,^{k+1}}{n-k} =
{n \choose k} (k+1)^{\frac{k-1}{2}},
$$
proving our assertion.\eop

\bigskip

Note that $C_n(A)=\det(A)$. Therefore, in the case $k=n$ the
Hadamard determinant theorem yields the following.

\begin{thm}
We have
$$
\theta_n(\ell_\infty,\ell_\infty) = \sqrt{n}\,^n.
$$
\end{thm}

\bigskip

In the cases $k=1$ and $k=n-1$ we have equality in (\ref{a11}) as
follows.

\begin{thm}
We have $\theta_1(\ell_\infty,\ell_\infty) = n.$
\end{thm}

{\it Proof.} In view of the inequality (\ref{a11}) all we have to
show is that there exists an $n \times n$ matrix $A$ satisfying
$||col_i(A)||_\infty = 1,\,\,i = 1,\ldots,n$ and such that
$||A||_\infty = n$. It is easy to check that the $n \times n$
matrix whose first row consists of $1$'s and all other entries
equal to $0$ is such a matrix.\eop

\bigskip

In order to establish the case $k=n-1$ we first make an
observation.

\begin{obs}
\rm For every positive integer $n$ there exists an $n \times n$
complex matrix $A$ satisfying $|a_{ij}|=1,\,\,i,j =1,\ldots,n$
and $AA^*=nI$. An example of such a matrix is the Vandermonde
matrix
$$
\left(\begin{array}{ccccc}
1 & 1 & 1 & \ldots & 1 \\
1 & \omega & \omega^2 & \ldots & \omega^{n-1} \\
1 & \omega^2 & \omega^4 & \ldots & \omega^{2n-2} \\
\vdots & \vdots & \vdots & & \vdots \\
1 & \omega^{n-1} & \omega^{2n-2} & \ldots & \omega^{(n-1)^2} \\
\end{array}\right),
$$
where $\omega=e^{\frac{2\pi i}{n}}$.
There are also the Hadamard matrices for those $n$'s for which
they exist.
\end{obs}

\begin{thm}
We have $\theta_{n-1}(\ell_\infty,\ell_\infty) = \sqrt{n}\,^n.$
\end{thm}

{\it Proof.} In view of the inequality (\ref{a11}) all we have to
show is that there exists an $n \times n$ matrix $A$ satisfying
$||col_i(A)||_\infty = 1,\,\,i = 1,\ldots,n$ and such that
$||C_{n-1}A||_\infty = \sqrt{n}\,^n$.
Let A be an $n \times n$ complex matrix satisfying
$|a_{ij}|=1,\,\,i,j =1,\ldots,n$ and $AA^*=nI$. Then
$A^{-1} = \frac{1}{n}A^*$, and so
$$
C_{n-1}(A) = (\det(A)DA^{-1}D)^T = \frac{det(A)}{n}D\overline AD
,
$$
where $D$ is the diagonal matrix with alternating $1$'s and
$-1$'s along the diagonal.
It now follows that $||C_{n-1}(A)||_\infty = \sqrt{n}\,^n$,
proving our claim.\eop

\bigskip

In order to consider some other combinations of norms, for a real
number $r$ we denote $[r]^+ = \max \{\,r\,,\,0\,\}$.

\begin{lem}
\label{l2}
Let $\mu$ be an absolute operator norm and let $p$ and $r$ be
positive integers. Then
$$
\theta_k(\mu,\ell_p) \leq
n^{\left[\frac{k}{r}-\frac{k}{p}\right]^+}\theta_k(\mu,\ell_r)
$$
\end{lem}

{\it Proof.} By \cite[p.26 \#16 and p.29 \#19]{HLP}, for every
vector $v$ in $\C^n$ we have
\beq
\label{hlp1}
||v||_r \leq n^{\left[\frac{1}{r}-\frac{1}{p}\right]^+}||v||_p.
\end{equation}

Our claim now follows from Theorem \ref{thm0} and from the fact
that (\ref{a0}) is sharp.\eop

\bigskip

From Theorems \ref{thml1}, \ref{tk2} and \ref{thm212} we obtain,
by Lemma \ref{l2}, the following corollary.

\begin{cor}
We have
\beq
\label{th12}
\theta_k(\ell_1,\ell_2) \leq \sqrt{n}\,^k ,
\end{equation}

\beq
\label{th1i}
\theta_k(\ell_1,\ell_\infty) \leq n^k ,
\end{equation}

\beq
\label{th21}
\theta_k(\ell_2,\ell_1) \leq
\left(\frac{n}{k}\right)^\frac{k}{2},
\end{equation}

\beq
\label{th2i}
\theta_k(\ell_2,\ell_\infty) \leq
\left(\frac{n}{k}\right)^\frac{k}{2}\sqrt{n}\,^k ,
\end{equation}

\beq
\label{thi1}
\theta_k(\ell_\infty,\ell_1) \leq
\left\{
\begin{array}{ll}
{n \choose k} (k+1)^{\frac{k-1}{2}}\,, & k<n \\
\\
\sqrt{n}\,^n\,, & k=n
\end{array}
\right. ,
\end{equation}

and

\beq
\label{thi2}
\theta_k(\ell_\infty,\ell_2) \leq
\left\{
\begin{array}{ll}
{n \choose k} (k+1)^{\frac{k-1}{2}}\,, & k<n \\
\\
\sqrt{n}\,^n\,, & k=n
\end{array}
\right. .
\end{equation}

\end{cor}

\bigskip

We conclude this section with two general remarks.

\begin{rem}
\rm Note that one can define $\theta_k (\mu,\nu)$ using rows
instead
of columns and obtain similar results where rows replace columns
all
over.
\end{rem}

\begin{rem}
\rm Denote by $adj(A)$ the
(classical) adjoint matrix of $A$, that is, the transposed matrix
of
cofactors. Note that the term {\it adjugate} is sometimes used
instead of {\it adjoint} to avoid confusion with the Hermitian
adjoint $A^*$. Since $adj(A)=DC_{n-1}(A^T)D$ where $D$ is the
diagonal matrix with alternating $1$'s and $-1$'s along the
diagonal, it follows that for absolute norms our results in the
case $k=n-1$ yield an upper bound on $\mu(adj(A))$. In
particular, our remark applies to the $\ell_1$, $\ell_2$ and
$\ell_\infty$ norms under discussion. Our results can also be
applied to the adjoint compounds found in \cite[Chapter 5]{Ait}.

\end{rem}

\section{Bounds on products of eigenvalues \label{lower}}

For an $n \times n$ complex matrix $A$ we denote by
$\lambda_1(A),\ldots,\lambda_n(A)$ the eigenvalues of $A$ ordered
in
a non-increasing order of their moduli. In this section we find
an
upper bound on the product $|\prod_{i=1}^k\lambda_i(A)|$ or,
equivalently, a lower bound on the product
$|\prod_{i=k+1}^n\lambda_i(A)|$.

\bigskip

Our results follow from the following corollary of Theorem
\ref{thm0}.
Here we denote by $row_i(A)$ the (vector in
$\C^n$ which is the) transpose of the $i$th row of $A$.

\begin{thm}
\label{thm1}
Let $\mu$ be an absolute operator norm on $\C^{mm}$ where $m={n
\choose k}$. Then
$$
\left|\prod_{i=1}^k\lambda_i(A)\right| \leq
\theta_k(\mu,\nu)\,
\min \left\{
\max_{\stackrel{\scriptstyle \alpha \subseteq
\{1,\ldots,n\}}{|\alpha|=k}} \prod_{i \in \alpha}\nu(col_i(A))
\ ,
\max_{\stackrel{\scriptstyle \alpha \subseteq
\{1,\ldots,n\}}{|\alpha|=k}} \prod_{i \in \alpha}\nu(row_i(A))
\right\}.
$$

\end{thm}

{\it Proof.}
As is well known, the spectral radius $\rho(C_k(A))$ of $C_k(A)$
satisfies
\beq
\label{ee1}
\rho(C_k(A))=\left| \prod_{i=1}^k\lambda_i(A)\right|.
\end{equation}

Also, we have $\rho(C_k(A))=\rho(C_k(A^T))$.
Since $\mu$ is an operator norm we have
\beq
\label{ee2}
\rho(C_k(A))\leq\mu(C_k(A)),
\end{equation}
Our claim follows from (\ref{ee1}), (\ref{ee2}) and
(\ref{a0}), where the latter is applied both to $A$ and
$A^T$.\eop

\bigskip

It follows from Theorem \ref{thm1} that in order to obtain a
better upper bound on $\left|\prod_{i=1}^k\lambda_i(A)\right|$ in
terms of the $\nu$ norms of the rows and columns of $A$, we
should pick us the $\mu$ norm that provides the lowest value of
$\theta_k(\mu,\nu)$. We now apply this approach to the results of
the previous section. The best upper bound in terms of the
$\ell_1$ norms of rows and columns that can be derived from
Theorem \ref{thml1} and from the inequalities (\ref{th21}) and
(\ref{thi1}) is

\begin{thm}
Let $A \in \C^{nn}$. Then
\beq
\label{in1}
\left|\prod_{i=1}^k\lambda_i(A)\right| \leq
\min \left\{
\max_{\stackrel{\scriptstyle \alpha \subseteq
\{1,\ldots,n\}}{|\alpha|=k}} \prod_{i \in \alpha}||col_i(A)||_1
\ ,
\max_{\stackrel{\scriptstyle \alpha \subseteq
\{1,\ldots,n\}}{|\alpha|=k}} \prod_{i \in \alpha}||row_i(A)||_1
\right\}.
\end{equation}
\end{thm}

\begin{rem}
\rm The upper bound on $|\prod_{i=1}^k\lambda_i(A)|$ given by
(\ref{in1}) is sharp since equality is attained for every
monomial matrix $A$.
\end{rem}

\begin{rem}
\rm
The inequality (\ref{in1}) was already proven in Theorem 1 of
\cite{Sch}, see also \cite[p.145, Theorem 1.7]{MM}, using
essentially the same techniques we do. The inequality (\ref{in1})
was also proven in Theorem 8 of \cite{SW}, and is weaker than
\cite[Formula (25)]{O}.
\end{rem}

\bigskip

The following theorem states the best upper bound on
$|\prod_{i=1}^k\lambda_i(A)|$ in terms of the $\ell_2$ norms of
rows and columns of $A$ that can be derived from Theorem
\ref{tk2} and from the inequalities (\ref{th12}) and
(\ref{thi2}).

\begin{thm}
Let $A \in \C^{nn}$. Then
\beq
\label{in2}
\left|\prod_{i=1}^k\lambda_i(A)\right| \leq
\left(\frac{n}{k}\right)^\frac{k}{2}\,
\min \left\{
\max_{\stackrel{\scriptstyle \alpha \subseteq
\{1,\ldots,n\}}{|\alpha|=k}} \prod_{i \in \alpha}||col_i(A)||_2
\ ,
\max_{\stackrel{\scriptstyle \alpha \subseteq
\{1,\ldots,n\}}{|\alpha|=k}} \prod_{i \in \alpha}||row_i(A)||_2
\right\}.
\end{equation}
\end{thm}

\begin{rem}
\rm In order to justify that (\ref{in2}) is indeed the best we
can derive from Theorem \ref{tk2} and from the inequalities
(\ref{th12}) and (\ref{thi2}), we have to show that for $k<n$ we
have
$$
\left(\frac{n}{k}\right)^\frac{k}{2} \leq
{n \choose k} (k+1)^{\frac{k-1}{2}}.
$$
This follows from a stronger inequality, see Remark \ref{add3}
below.
\end{rem}

\begin{rem}
\rm Note that by Theorem \ref{tk2} the inequality (\ref{in2}) is
strict whenever $A$ is a nonsingular matrix.
\end{rem}

\begin{rem}
\rm The upper bound on $|\prod_{i=1}^k\lambda_i(A)|$ given by
(\ref{in2}) is sharp. It is easy to check that, as in the proof
of the equality case in Theorem \ref{tk2}, equality holds for a
positive semidefinite $n \times n$ matrix $A$ such that the
positive semidefinite matrix $B = A^2$ has diagonal elements all
equal to $1$ and where the eigenvalues of $B$ are $\frac{n}{k}$
with multiplicity $k$ and $0$ with multiplicity $n-k$.
\end{rem}

\begin{rem}
\rm Inequality (\ref{in2}) is a generalization of Hadamard
determinant theorem, which is the special case of (\ref{in2})
where $k=n$.
\end{rem}

\begin{rem}
\rm The special case of (\ref{in2}) where $k=n-1$, that is, the
inequality
$$
\left|\prod_{i=1}^{n-1}\lambda_i(A)\right| \leq
\left(1 + \frac{1}{n-1}\right)^{\frac{n-1}{2}}
\min \left\{
\max_k\prod_{\stackrel{\scriptstyle j=1}{j\neq
k}}^n||row_j(A)||_2\,,\,
\max_k\prod_{\stackrel{\scriptstyle j=1}{j\neq
k}}^n||col_j(A)||_2
\right\} ,
$$
follows from \cite[Theorem 1]{HP}.
\end{rem}

\begin{rem}
\rm
Another result in Theorem 8 of \cite{SW} could be stated as
$$
\left|\prod_{i=1}^{n-1}\lambda_i(A)\right| \leq
n\,\max_k\prod_{\stackrel{\scriptstyle j=1}{j\neq
k}}^n||row_j(A)||_2.
$$
This inequality follows from the special case of our inequality
(\ref{in2}) where $k=n-1$. It is, in fact, weaker than our result
since
$$
\left(1 + \frac{1}{n-1}\right)^{\frac{n-1}{2}} < \sqrt{e} < n
\,\,\,\,\,\,\,\,\,\,\,\,\,\,\,(n>1).
$$

\end{rem}

\begin{rem}
\label{r312}
\rm The upper bounds on $|\prod_{i=1}^k\lambda_i(A)|$ given by
(\ref{in1}) and by (\ref{in2}) are not comparable. The bound
given by (\ref{in1}) is
better, for example, in the case of a monomial matrix $A$, since
in such a
case the $\ell_1$ norm and the $\ell_2$ norm of the rows (and
columns) of $A$ are the same. On the other hand, if $A$ is an $n
\times n$ complex matrix satisfying $|a_{ij}|=1,\,\,i,j =
1,\ldots,n$ and $AA^*=nI$ then the
$\ell_1$ norm of any row and column  of $A$ is equal to $n$,
while
the $\ell_2$ norm of any row and column of $A$ is equal to
$\sqrt{n}$.
Therefore, the left hand side of (\ref{in1}) becomes $n^k$ while
the left hand side of (\ref{in2}) becomes
$\frac{n^k}{k^{\frac{k}{2}}}$, which is a better upper bound.
\end{rem}

\bigskip

The following theorem states the best upper bound on
$|\prod_{i=1}^k\lambda_i(A)|$ in terms of the $\ell_\infty$ norms
of rows and columns of $A$ that can be derived from Theorem
\ref{thm212} and from the inequalities (\ref{th1i}) and
(\ref{th2i}).

\begin{thm}
Let $A \in \C^{nn}$. Then
\beq
\label{in3}
\left|\prod_{i=1}^k\lambda_i(A)\right| \leq
\end{equation}
$$
\left(\frac{n}{k}\right)^\frac{k}{2}\sqrt{n}\,^k\,
\min \left\{
\max_{\stackrel{\scriptstyle \alpha \subseteq
\{1,\ldots,n\}}{|\alpha|=k}} \prod_{i \in
\alpha}||col_i(A)||_\infty
\ ,
\max_{\stackrel{\scriptstyle \alpha \subseteq
\{1,\ldots,n\}}{|\alpha|=k}} \prod_{i \in
\alpha}||row_i(A)||_\infty
\right\}.
$$
\end{thm}

\begin{rem}
\label{add3}
\rm In order to justify that (\ref{in3}) is indeed the best we
can derive from Theorem \ref{thm212} and from the inequalities
(\ref{th1i}) and (\ref{th2i}), we have to show that for $k<n$ we
have
$$
\left(\frac{n}{k}\right)^\frac{k}{2}\sqrt{n}\,^k \leq
{n \choose k} (k+1)^{\frac{k-1}{2}}
$$
or, equivalently,
\beq
\label{aa1}
h(k,n)=
\frac{\left(\frac{n}{k}\right)^\frac{k}{2}\sqrt{n}\,^k}
{{n \choose k} (k+1)^{\frac{k-1}{2}}}=
\frac{n^k}{\sqrt{k}\,^k \sqrt{k+1}\,^{k-1}{n \choose k}} \leq 1.
\end{equation}
Note that
\beq
\label{aa2}
h(1,n)=1,\,\,\,\forall n.
\end{equation}
Since $k \leq n-1$ we have $(2n+1)k < 2n^2$, which is equivalent
to
$$
\frac{k+2}{k} > \frac{(n+1)^2}{n^2},
$$
or
\beq
\label{aa3}
\frac{n\,\sqrt{k+2}}{\sqrt{k}\,(n+1)} > 1.
\end{equation}
It is easy to check that
\beq
\label{aa4}
\frac{h(k,n)}{h(k+1,n+1)} =
\left(\frac{n\,\sqrt{k+2}}{\sqrt{k}\,(n+1)}\right)^k.
\end{equation}
It now follows from (\ref{aa2}), (\ref{aa3}) and (\ref{aa4}) that
for every $k$ and $n$, $k \leq n$, we have $h(k,n) \leq 1$,
proving (\ref{aa1}).
\end{rem}

\begin{rem}
\rm By (\ref{hlp1}), the upper bound on
$|\prod_{i=1}^k\lambda_i(A)|$ given by (\ref{in3}) follows from
the one given by (\ref{in2}). The bounds given by (\ref{in3}) and
by (\ref{in1}) are not comparable. The bound given by (\ref{in1})
is better, for example, in the case of a monomial matrix $A$,
since in such a case the $\ell_1$ norm and the $\ell_\infty$ norm
of the rows (and columns) of $A$ are the same. On the other hand,
if $A$ is an $n \times n$ complex matrix satisfying
$|a_{ij}|=1,\,\,i,j = 1,\ldots,n$ and $AA^*=nI$ then the $\ell_1$
norm of any row and column of $A$ is equal to $n$, while the
$\ell_\infty$ norm of any row and column of $A$ is equal to $1$.
Therefore, the left hand side of (\ref{in1}) becomes $n^k$ while
the left hand side of (\ref{in3}) becomes
$\frac{n^k}{k^{\frac{k}{2}}}$, which is a better upper bound.
\end{rem}

Our final remark refers to products of smallest (moduli of)
eigenvalues of a given $n \times n$ matrix $A$.

\begin{rem}
\rm Since $\prod_{i=k+1}^n\lambda_i(A) =
\frac{\det(A)}{\prod_{i=1}^k\lambda_i(A)}$, it follows that all
the results of this section on upper bounds on the products
$|\prod_{i=1}^k\lambda_i(A)|$ of $k$ largest eigenvalues of $A$
yield, whenever $A$ is nonsingular, lower bounds on the products
$|\prod_{i=k+1}^n\lambda_i(A)|$ of $n-k$ smallest eigenvalues of
$A$.
\end{rem}

\end{document}